\documentclass[10pt]{amsart}
\usepackage{amsmath,amssymb,latexsym,hyperref,url,graphicx}

\newtheorem*{theorem}{Theorem}

\newcommand{\Z}{\mathbb{Z}}

\begin{document}

\title{Non-regularity of $\lfloor \alpha + \log_k n \rfloor$}
\author{Eric S. Rowland}
\address{
	Mathematics Department \\
	Tulane University \\
	New Orleans, LA 70118, USA
}
\date{October 11, 2009}

\thanks{Thanks are due to the referee for a careful reading and corrections.}

\begin{abstract}
This paper presents a new proof that if $k^\alpha$ is irrational then the sequence $\{\lfloor \alpha + \log_k n \rfloor\}_{n \geq 1}$ is not $k$-regular.  Unlike previous proofs, the methods used do not rely on automata or language theoretic concepts.  The paper also proves the stronger statement that if $k^\alpha$ is irrational then the generating function in $k$ non-commuting variables associated with $\{\lfloor \alpha + \log_k n \rfloor\}_{n \geq 1}$ is not algebraic.
\end{abstract}

\maketitle
\markboth{Eric Rowland}{Non-regularity of $\lfloor \alpha + \log_k n \rfloor$}

Fix an integer $k \geq 2$.  A sequence $\{a(n)\}_{n \geq 0}$ is \emph{$k$-regular} if the $\Z$-module generated by the subsequences $\{a(k^e n + i)\}_{n \geq 0}$ for $e \geq 0$ and $0 \leq i < k^e$ is finitely generated.  Regular sequences were introduced by Allouche and Shallit \cite{Allouche-Shallit} and have several nice characterizations, including the following characterization as rational power series in non-commuting variables $x_0, x_1, \dots, x_{k-1}$.  If $n = n_l \cdots n_1 n_0$ is the standard base-$k$ representation of $n$, then let $\tau(n) = x_{n_0} x_{n_1} \cdots x_{n_l}$.  The sequence $\{a(n)\}_{n \geq 0}$ is $k$-regular if and only if the power series $\sum_{n \geq 0} a(n) \tau(n)$ is rational.  In this sense, regular sequences are analogous to constant-recursive sequences (sequences that satisfy linear recurrence equations with constant coefficients), the set of which coincides with the set of sequences whose generating functions in a single variable are rational.

The sequence $\{\lfloor \log_2 (n + 1) \rfloor\}_{n \geq 0}$ is an example of a $2$-regular sequence, and the associated power series in non-commuting variables $x_0$ and $x_1$ is
\begin{multline*}
	f(x_0, x_1) = \sum_{n \geq 0} \lfloor \log_2 (n+1) \rfloor \tau(n) \\
	= x_1 + x_0 x_1 + 2 x_1 x_1 + 2 x_0 x_0 x_1 + 2 x_1 x_0 x_1 + 2 x_0 x_1 x_1 + 3 x_1 x_1 x_1 + \cdots.
\end{multline*}
The rational expression for this series is somewhat large; however its commutative projection is quite manageable:
\[
	\frac{x_1 \left(1 - x_0 - x_1 + x_0^2 + x_0 x_1\right)}{\left(1 - x_1\right) \left(1 - x_0 - x_1\right)^2}.
\]

Allouche and Shallit \cite[open problem~16.10]{Automatic Sequences} asked whether the sequence $\{\lfloor \frac{1}{2} + \log_2 (n + 1) \rfloor\}_{n \geq 0}$ is $2$-regular.  Bell \cite{Bell} and later Moshe \cite[Theorem~4]{Moshe} gave proofs that this sequence is not $2$-regular.  Moreover, they proved the following.

\begin{theorem}
Let $k \geq 2$ be an integer and $\alpha$ be a real number.  The sequence $\{\lfloor \alpha + \log_k (n+1) \rfloor\}_{n \geq 0}$ is $k$-regular if and only if $k^\alpha$ is rational.
\end{theorem}

In this paper we prove the following theorem, which is a slightly weaker statement than the previous theorem but still establishes that if $k^\alpha$ is irrational then $\{\lfloor \alpha + \log_k (n+1) \rfloor\}_{n \geq 0}$ is not $k$-regular.  Let $|\tau(n)|$ be the length of the word $\tau(n)$, i.e., $|\tau(0)| = 0$ and $|\tau(n)| = \lfloor \log_k n \rfloor + 1$ for $n \geq 1$.

\begin{theorem}
Let $k \geq 2$ be an integer and $\alpha$ be a real number.  The series $f(x) = \sum_{n \geq 0} \lfloor \alpha + \log_k (n+1) \rfloor x^{|\tau(n)|}$ is rational if and only if $k^\alpha$ is rational.
\end{theorem}

The proof given here is similar to Moshe's but does not require the notion of a regular language.  Note that, given the associated power series
\[
	f(x_0, x_1, \dots, x_{k-1}) = \sum_{n \geq 0} \left\lfloor \alpha + \log_k (n+1) \right\rfloor \tau(n),
\]
the series in the theorem is the power series $f(x) = f(x, x, \dots, x)$ in one variable obtained by setting $x_0 = x_1 = \cdots = x_{k-1} = x$.  Therefore non-rationality of $f(x)$ implies non-regularity of $\{\lfloor \alpha + \log_k (n+1) \rfloor\}_{n \geq 0}$.

To get a sense of computing $f(x)$ in the proof of the theorem, first we examine the case where $k = 2$ and $\alpha = \frac{1}{2}$.  The power series in this case is
\begin{multline*}
	f(x_0, x_1) = \sum_{n \geq 0} \left\lfloor \frac{1}{2} + \log_2 (n+1) \right\rfloor \tau(n) \\
	= x_1 + 2 x_0 x_1 + 2 x_1 x_1 + 2 x_0 x_0 x_1 + 3 x_1 x_0 x_1 + 3 x_0 x_1 x_1 + 3 x_1 x_1 x_1 + \cdots,
\end{multline*}
and
\begin{align*}
	f(x) &= \sum_{n \geq 0} \left\lfloor \frac{1}{2} + \log_2 (n+1) \right\rfloor x^{|\tau(n)|} \\
	&= x + 2 x^2 + 2 x^2 + 2 x^3 + 3 x^3 + 3 x^3 + 3 x^3 + 3 x^4 + 3 x^4 + 3 x^4 + 4 x^4 + \cdots \\
	&= x + 4 x^2 + 11 x^3 + 29 x^4 + 74 x^5 + 179 x^6 + 422 x^7 + 971 x^8 + 2198 x^9 + \cdots \\
	&= \sum_{m \geq 0} b(m) x^m.
\end{align*}
To write $b(m)$ in closed form, we observe how the first few terms of $\{\lfloor \frac{1}{2} + \log_2 (n + 1) \rfloor\}_{n \geq 0}$ gather by exponent:
\[
	0 \mspace{1mu} 1 \mspace{1mu} 22 \mspace{-1mu} \underbrace{2333}_{x^3} \underbrace{33344444}_{x^4} \underbrace{4444445555555555}_{x^5} \underbrace{55555555555556666666666666666666}_{x^6} \cdots
\]
Since the length of $n$ in binary is $|\tau(n)| = 1 + \lfloor \log_2 n \rfloor$ for $n \geq 1$, the difference $|\tau(n)| - \lfloor \frac{1}{2} + \log_2 (n+1) \rfloor$ between exponent and coefficient in each term of the first sum above is either $1$ or $0$.  In other words, the only terms that contribute to $b(m) x^m$ are of the form $(m-1) x^m$ and $m x^m$, so for some sequence $\{c(m)\}_{m \geq 1}$ we have
\[
	b(m) = (m-1) \left(c(m) - 2^{m-1}\right) + m \left(2^m - c(m)\right)
\]
for $m \geq 1$.  In fact $c(m)$ is the smallest value of $n$ for which $\frac{1}{2} + \log_2 (n+1) \geq m$, so $c(m) = \lfloor 2^{m - \frac{1}{2}} \rfloor$ and $b(m) = (m+1) 2^{m-1} - \lfloor 2^{m - \frac{1}{2}} \rfloor$ for $m \geq 1$.  Therefore
\[
	f(x) = \frac{1}{2 (1 - 2 x)^2} - \frac{1}{2} - \sum_{m \geq 0} \left\lfloor 2^{m - \frac{1}{2}} \right\rfloor x^m,
\]
where the term $-1/2$ is needed because $b(0) = 0$.

We carry out the preceding computation more generally to prove the theorem.

\begin{proof}
Let $\text{frac}(\alpha) = \alpha - \lfloor\alpha\rfloor$ denote the fractional part of $\alpha$.  Then
\begin{align*}
	f(x) &= \sum_{n \geq 0} \left\lfloor \alpha + \log_k (n+1) \right\rfloor x^{|\tau(n)|} \\
	&= \left\lfloor \alpha + \log_k 1 \right\rfloor + \sum_{m \geq 1} \sum_{i = k^{m-1}}^{k^m - 1} \left\lfloor \alpha + \log_k (i+1) \right\rfloor x^m \\
	&= \left\lfloor \alpha \right\rfloor + \sum_{m \geq 1} \left( \sum_{i = k^{m-1}}^{\left\lceil k^{m - \text{frac}(\alpha)} \right\rceil - 2} \left\lfloor \alpha + \log_k (i+1) \right\rfloor + \sum_{i = \left\lceil k^{m - \text{frac}(\alpha)} \right\rceil - 1}^{k^m-1} \left\lfloor \alpha + \log_k (i+1)\right\rfloor \right) x^m.
\end{align*}
Since
\[
	\left\lfloor \alpha + \log_k (i+1) \right\rfloor =
	\begin{cases}
		\lfloor \alpha \rfloor + m - 1	& \text{if $k^{m-1} + 1 \leq i + 1 \leq \left\lceil k^{m - \text{frac}(\alpha)} \right\rceil - 1$} \\
		\lfloor \alpha \rfloor + m	& \text{if $\left\lceil k^{m - \text{frac}(\alpha)} \right\rceil \leq i + 1 \leq k^m$},
	\end{cases}
\]
we have
\begin{align*}
	f(x) &= \left\lfloor \alpha \right\rfloor + \sum_{m \geq 1} \left( k^{m-1} \left((k-1) (m + \lfloor \alpha \rfloor) + 1\right) + 1 - \left\lceil k^{m-\text{frac}(\alpha)} \right\rceil \right) x^m \\
	&= \frac{(1 - x) (k x + \lfloor \alpha \rfloor (1 - k x))}{(1 - k x)^2} + \frac{x}{1-x} + \sum_{m \geq 1} \left\lfloor -k^{m - \text{frac}(\alpha)} \right\rfloor x^m.
\end{align*}
The series $f(x)$ is therefore rational if and only if
\begin{align*}
	g(x) &= -\left\lfloor -k^{1 - \text{frac}(\alpha)} \right\rfloor + \left(\frac{1}{x} - k\right) \sum_{m \geq 1} \left\lfloor -k^{m - \text{frac}(\alpha)} \right\rfloor x^m \\
	&= \sum_{m \geq 1} \left( \left\lfloor -k^{m + 1 - \text{frac}(\alpha)} \right\rfloor - k \left\lfloor -k^{m - \text{frac}(\alpha)} \right\rfloor \right) x^m
\end{align*}
is rational.  The expression $\lfloor k^m y \rfloor - k \lfloor k^{m-1} y \rfloor$ is the $(-m)$th base-$k$ digit of $y$, so the coefficients of $g(x)$ are the base-$k$ digits of $\text{frac}(-k^{1-\text{frac}(\alpha)})$, which is rational precisely when $k^\alpha$ is rational.

If $k^\alpha$ is rational, then the coefficients of $g(x)$ are eventually periodic, so $g(x)$ and hence $f(x)$ is rational.  If $k^\alpha$ is irrational, then $g(x)$ is not rational, since in particular $g(\frac{1}{k}) = \text{frac}(-k^{1-\text{frac}(\alpha)})$ is irrational; therefore $f(x)$ is not rational.
\end{proof}

In fact we may show something stronger:  Not only does $f(x_0, x_1, \dots, x_{k-1})$ fail to be rational when $k^\alpha$ is irrational, but it fails to be algebraic.  Bell, Gerhold, Klazar, and Luca \cite[Proposition~13]{Bell-Gerhold-Klazar-Luca} prove that if a polynomial-recursive sequence (a sequence satisfying a linear recurrence equation with polynomial coefficients) has only finitely many distinct values, then it is eventually periodic.  It follows that the coefficient sequence of $g(x)$ is not polynomial-recursive, hence $g(x)$ is not algebraic, and $f(x, x, \dots, x)$ is not algebraic.

\end{document}